\documentclass{article}
\usepackage{amsmath,amssymb, pxfonts,graphicx,color}
\numberwithin{equation}{section}
\newtheorem{thm}{Theorem}[section]
\newtheorem{lem}[thm]{Lemma}
\newtheorem{remark}[thm]{Remark}

\newcommand{\re}{\mathbb{R}}
\newcommand{\E}{\mathbb{E}}

\usepackage{hyperref}

\title{Efficient discretisation of stochastic differential equations}
\author{Masaaki Fukasawa\thanks{Graduate School of Engineering Science,
and Center for Mathematical modeling and Data Science, Osaka University,
1-3 Machikaneyama,
Toyonaka, Osaka, Japan. email:
\href{mailto:fukasawa@sigmath.es.osaka-u.ac.jp}\nolinkurl{fukasawa@sigmath.es.osaka-u.ac.jp},
web: \url{www.sigmath.es.osaka-u.ac.jp/\~fukasawa/}. Masaaki Fukasawa's 
work was supported by JSPS KAKENHI Grant Number 17K05297.} 
 \and Jan Ob\l\'oj\thanks{Mathematical Institute, University of Oxford, ROQ, Woodstock Road, Oxford OX2 6GG, UK. email: \href{mailto:jan.obloj@maths.ox.ac.uk}\nolinkurl{jan.obloj@maths.ox.ac.uk}, web: \url{www.maths.ox.ac.uk/people/jan.obloj}.
 Jan Ob\l\'oj gratefully acknowledges funding received from the European Research Council under the European Union's Seventh Framework Programme (FP7/2007-2013) / ERC grant agreement no. 335421 and is also thankful to the Oxford-Man Institute of Quantitative Finance and St John's College in Oxford for their financial support.
}}

\date{}

\begin{document}
\maketitle
\begin{abstract}
The aim of this study is to find a generic method for generating a
path of the solution of a given stochastic differential equation  which is
more efficient than the standard Euler--Maruyama scheme with 
Gaussian increments. 
First we characterize the asymptotic distribution of pathwise error in 
 the Euler--Maruyama scheme with a general partition of time interval and then,
show that the error is reduced by a factor 
$(d+2)/d$ when using a partition associated with the hitting times of sphere for the
driving $d$--dimensional Brownian motion.
This reduction ratio is the best possible in a symmetric class of partitions.
Next we show that a reduction which is close to the best possible is achieved by
using the hitting time of a  moving sphere that is easier to implement.
\end{abstract}

\section{Introduction}
Various stochastic phenomena have been modeled  in terms of
the solution $X = (X^1,\dots,X^p)$ of a stochastic differential equation (SDE)
\begin{equation} \label{sde}
\mathrm{d}X^i(s) = \sum_{j=0}^d f^i_j(X(s))\mathrm{d}W^j(s)
\end{equation}
on a domain $\mathbb{D} \subset \mathbb{R}^p$,
where $W = (W^1,\dots,W^d)$ is a $d$-dimensional standard Brownian
motion, $W^0(s) = s$ is the time coordinate, and 
$f = \{f^i_j; 1 \leq i \leq p, 0 \leq j \leq d\} : \mathbb{D} \to \mathbb{R}^p \otimes \mathbb{R}^{d+1}$, 
is a continuously differentiable function.
The Monte Carlo simulation is a powerful and very popular approach to study such a
stochastic model. The standard method for generating a path which
follows the SDE is the Euler--Maruyama scheme that 
constructs an approximating sequence of processes
$X^n = (X^{n,1},\dots,X^{n,p})$ 
as
\begin{equation}\label{em}
X^n(0) = X(0), \ \ 
X^{n,i}(s) = X^{n,i}({\pi^n_m}) + \sum_{j=0}^{d}f^i_j(X^{n}({\pi^n_m}))
(W^j(s)-W^j({\pi^n_m})) 
\end{equation} 
for $s\in (\pi^n_m, \pi^n_{m+1}]$,
where $\pi^n = \{\pi^n_j\}$ is an increasing sequence of stopping
times, usually chosen to be a deterministic sequence such as
$\pi^n_m = m/n$.
The variable $n$ controls the computational effort 
of this construction. We naturally expect that $X^n \to X$ in some sense 
as $n \to \infty$. 
The attractive features of the Euler--Maruyama
scheme include its validity under degenerate 
diffusion coefficients with
mild regularity,
intuitive construction and easy implementation.
See Kloeden and Platen~\cite{KP} for some elementary properties of
this and other related methods.

Kurtz and Protter~\cite{KP2} studied the limit of the
approximation error process
\begin{equation}\label{ep}
U^n = (U^{n,1},\dots, U^{n,p}), \ \ 
 U^{n,i}(s) = \sqrt{n}(X^{n,i}(s)-X^i(s))
 \end{equation}
 as $n\to \infty$. They showed that
if the sequence of $(d+1)^2$-dimensional processes 
$Z^n = \{Z^{n,l,j}; 0\leq l \leq d, 0\leq j \leq d\}$ defined by
\begin{equation}\label{defz}
Z^{n,l,j}(t)
= \sqrt{n}\sum_{m=0}^{\infty} \int_{t \wedge \pi^n_m}^{t \wedge \pi^n_{m+1}}
(W^l(s) - W^l({\pi^n_m}))\mathrm{d}W^j(s)
\end{equation}
is ``good''  and converging to a semimartingale $Z = \{Z^{l,j}\}$ in law, 
then $U^{n}$ converges in law
to $U = (U^1,\dots,U^p)$, the solution of
\begin{equation}\label{SDEU}
\mathrm{d}U^i(s)
 =\sum_{j,k}\partial_kf^i_j(X(s))U^k(s) \mathrm{d}W^j(s) - \sum_{j,k,l}
\partial_kf^i_j(X(s))f^k_l(X(s))\mathrm{d}Z^{l,j}(s).
\end{equation}
Since this SDE for $U$ is affine, we may write $U$ in a
more explicit form. 
In particular if $\langle Z^{l,j}, W^i \rangle = 0$ for all $l,j,i$, 
then\footnote{ Theorem~56 of Chapter V,
Protter~\cite{Protter}.
Here, $U$ and $F$ are interpreted as column vectors.},
\begin{equation}\label{Uint}
U(t) = Y(t) \int_0^t Y(s)^{-1} \mathrm{d}F(s), \ \ 
\mathrm{d}F^i(s) =   - 
\sum_{j,k,l}\partial_kf^i_j(X(s))f^k_l(X(s))\mathrm{d}Z^{l,j}(s),
\end{equation}
where $Y = \{Y^{a,b}; 1 \leq a, b \leq p\}$ is the solution of
\begin{equation} \label{SDEY}
\mathrm{d}Y^{a,b}(s) = \sum_{j,k}
\partial_k f^a_j(X(s))Y^{k,b}(s) \mathrm{d}W^j(s),
\ \ Y^{a,b}(0) = \delta^{a,b},
\end{equation}
and $\delta^{a,b}=\mathbf{1}_{a=b}$ is Kronecker's delta. For example in the case of  $p=1$,
\begin{equation*}
U(t) =  - Y(t) \sum_{j,l=0}^d \int_0^t Y(s)^{-1} 
f_j^\prime (X(s)) f_l(X(s))\mathrm{d}Z^{l,j}(s),
\end{equation*}
and
\begin{equation*}
Y(t) = \exp\left\{ 
\int_0^t f_0^\prime (X(s))\mathrm{d}s  + 
\sum_{j=1}^d \left\{\int_0^t f_j^\prime (X(s))\mathrm{d}W^j(s) - 
\frac{1}{2}\int_0^t f_j^\prime(X(s))^2\mathrm{d}s\right\}\right\}.
\end{equation*}
It is important to note that $Y$ does not depend on $\pi^n$.
In consequence, the limit distribution $U$ of $U^n$ 
is determined by the limit distribution  $Z$ of $Z^n$ in a linear manner.
In the equidistant case $\pi^n_m = m/n$, as shown by
Kurtz and Protter~\cite{KP2} and  Jacod and Protter~\cite{JP},
we have 
$Z^{0,j} =Z^{j,0} = 0$ for all $ 0 \leq j \leq d$ and
$$\sqrt{2} Z^{l,j} , \ \ 1 \leq l,j \leq d$$ is a $d^2$-dimensional 
Brownian motion independent of $W$.
In particular, we have
 $\langle Z^{l,j}, W^i \rangle = 0$ for all $l,j,i$.
and that the distribution of $U(t)$ is conditionally 
Gaussian. 

In this paper we consider more general sequences of partitions $\{\pi^n\}$. Our first main result states a central limit theorem for the error process and provides a convenient characterization of the resulting process $Z$. This allows us to study the efficiency of different sequences of partitions. Our second main results establishes a uniform lower bound on the expected asymptotic error $U$ and shows that the bound is attained by partitions associated 
with the hitting times of a sphere for $W$. Note that in such a scheme both the time steps 
$\pi_{m}-\pi_{m-1}$ and the $W$ increments, $W({\pi_{m}})-W({\pi_{m-1}})$ are random. The later are uniformly distributed on a sphere and in particular in the one-dimensional case take just two values. 
In contrast, in the equidistant case the time steps are deterministic (and equal) and all the randomness is coming from the increments of $W$. 
Newton~\cite{Newton} and Fukasawa~\cite{F} studied the hitting time scheme in the one-dimensional case. Cambanis and Hu~\cite{CH} and 
M$\ddot{\text{u}}$ller-Gronbach~\cite{MG, MG2} gave optimality results among
irregular and adaptive schemes which still use (conditionally) 
Gaussian increments.
{ Our framework of discretisation admits multi-dimensional 
non-Gaussian (not even
conditionally) increments
$W(\pi_m)-W(\pi_{m-1})$ and therefore is not covered by
 these preceding studies.}
Our result is closely related to a recent work by Landon~\cite{Landon}, where 
a sequence of hitting times of ellipsoids is derived as an asymptotically optimal scheme. In this paper, we restrict schemes to be symmetric in a certain sense because, among other reasons, there is unlikely to be a realistic computational algorithm to implement asymmetric schemes. Our framework therefore excludes hitting times of ellipsoids (except spheres).

Our theoretical analysis of efficiency described above does not take into account the complexity of simulating random variables with a given distribution, which may be challenging for hitting times of spheres in higher dimensions. We argue that in practice one should look for a scheme where both the time steps and the spatial increments are random and all are easily generated together. We achieve this adapting the moving sphere approach in Deaconu and Herrmann~\cite{D}. 
This scheme is easier to implement and enjoys a better accuracy  than the
standard Gaussian scheme. In fact, it modifies the standard one only by replacing
\begin{equation*}
 \pi_{m}-\pi_{m-1} = \Delta, \ \
  W({\pi_m})-W({\pi_{m-1}}) = \sqrt{\Delta}N_m
\end{equation*}
with
\begin{equation*}
 \pi_{m}-\pi_{m-1} = \Delta^\prime e^{-Z_m}, \ \
  W({\pi_m})-W({\pi_{m-1}}) = \sqrt{\Delta^\prime d Z_me^{-Z_m}}\frac{N_m}{|N_m|},
\end{equation*}
where { $\Delta$ and $\Delta^\prime$ are constants to control
computational efforts,}
\begin{equation*}
  Z_m = \frac{|N_m|^2 + 2E_m}{d}
\end{equation*}
and $E_m \sim \mathrm{Exp}(1)$ and $N_m \sim \mathcal{N}(0, I_d)$ are
independent iid sequences. 
It improves the accuracy of the Monte Carlo
 simulation even after taking into account
a slight increase of computational time due to the one additional
 generation of exponential random variable and the calculation of the exponential function each step. 
It also has the further advantage that both the time and $W$ increments are bounded. 
This enables us to control the size of each increments of $X^n$ and deal with SDE on a bounded domain, see Milstein and Tretyakov~\cite{MT}, or devise efficient pricing of path dependent options e.g.\ barrier options.

This paper is organised as follows. 
In Section~2 we present a central limit theorem for
discretisation error. 
In Section~3 we study several examples of schemes 
and discuss their effectiveness.
In particular,  we show one of them to be attractive in terms of both  error magnitude
and computational costs.

\section{Central limit theorem for the error process}
Here we present a central limit theorem for the asymptotic error process $U^n = \sqrt{n}(X^n-X)$.
Let $(\Omega,\mathcal{F},\mathbb{P})$ be a probability space.
We suppose that the SDE (\ref{sde}) admits a unique strong solution $X$ which 
does not explode and remains in a given open connected domain $\mathbb{D}\subset \mathbb{R}^p$. We further assume there exists a sequence of compact sets 
$\mathbb{K}_m$ with each $\mathbb{K}_m$ being a subset of the interior of $\mathbb{K}_{m+1}$ such that
$\cup_{m=1}^{\infty}\mathbb{K}_m = \mathbb{D}$ and 
$\tau_m \to \infty$ as $m\to \infty$, where
\begin{equation*}
 \tau_m = \inf \{t > 0; X(t) \not\in \mathbb{K}_m\}.
\end{equation*}
Denote by $\{\mathcal{F}_s\}_{s \geq 0}$ the augmentation of the natural 
filtration generated by $W$.
A partition $\pi = \{\pi_m\}_{m \geq 0}$ is
 a sequence of increasing stopping times with $\pi_0 = 0$ and $\lim_{m \to \infty} \pi_m = \infty$. 
For  a partition $\pi$ and a finite stopping time $\tau$, put
\begin{equation*}
\begin{split}
& \Delta_m \pi = \pi_{m}-\pi_{m-1},  \\ 
&\Delta^\pi_mW = W({\pi_{m}})-W({\pi_{m-1}}), \\
& \mathcal{F}^\pi_m = \mathcal{F}_{\pi_m}, \\ 
& N^\pi_{\tau} = \min\{ m \geq 0; \pi_m \geq \tau\}.
\end{split}
\end{equation*} 
Note that $N^\pi_{\tau}$ is the number of discretisation 
steps required to generate a path up to a finite stopping time $\tau$.
In this section we do not take the computational difficulty to 
generate $\Delta_m \pi$
and $\Delta^\pi_mW$ into account. 
Hence we take $N^\pi_{\tau}$ as a measure of computational effort
associated with the partition $\pi$.
Notice that $N^\pi_\tau$ is a stopping time with respect to the discrete
filtration $\{\mathcal{F}^\pi_m\}_{m \geq 0}$.

For a given sequence of partitions  $\pi^n$,
we denote by $\Delta^n_m W$,  $N^n_\tau$ and $\mathcal{F}^n_m$ 
the corresponding quantities $\Delta^{\pi^n}_m W$, $N^{\pi^n}_\tau$ and
$\mathcal{F}_{\pi_m^n}$ for brevity. We also let $\pi^n(s) = \pi^n_m$ for $s \in [\pi^n_m,\pi^n_{m+1})$.
For example if $\{\pi^n\}$ is the equidistant scheme, that is,  $\pi^n_m = m/n$, then 
we have $N^n_\tau = \lceil n\tau \rceil \sim n\tau$ and 
$\pi^n(s) =\lfloor ns \rfloor/n \sim s$ as $n\to \infty$.
To see what happens if $\{\pi^n\}$ is stochastic, 
let us first consider an adaptive scheme.
Let $G$ be a positive continuous adapted process and define $\pi^n_m$ by 
\begin{equation} \label{pred}
\pi^n_0 = 0, \ \ 
\pi^n_{m + 1} = \pi^n_m + \frac{1}{n G({\pi^n_m})}. 
\end{equation}
Then
\begin{equation*}
\frac{N^n_\tau}{n} = \sum_{m=1}^{N^n_\tau}
G(\pi^n_{m-1}) \Delta_m \pi^n  \sim \int_0^\tau
G(s)\mathrm{d}s
\end{equation*}
as $n \to \infty$. 

More generally, if { there exists a locally integrable process $G$ such that
\begin{equation}\label{sumsq}
\sum_{m=1}^{N^n_\tau}\mathbb{E}\left[|\Delta_{m}\pi^n|^2 |
\mathcal{F}^n_{m-1}\right] \to 0, \ \ 
\sup_{ 0 \leq s \leq \tau }|G^n(\pi^n(s))- G(s)| \to 0
\end{equation}
in probability, where
\begin{equation*} 
G^n(\pi^n_m) =  \frac{1}{n
\mathbb{E}\left[\Delta_{m+1}\pi^n | \mathcal{F}^n_{m}\right]}, \ \ m=0,
1, \dots,
\end{equation*}
 then we have
\begin{equation} \label{effort}
\frac{N^n_\tau}{n} = \sum_{m=1}^{N^n_\tau} 
G^n(\pi^n_{m-1})\mathbb{E}\left[\Delta_{m}\pi^n | \mathcal{F}^n_{m-1}\right]
\to  \int_0^\tau
G(s)\mathrm{d}s
\end{equation}}
in probability as $n \to \infty$ by a simple application of the Lenglart
inequality (see e.g.\ Lemma~A.2 in Fukasawa~\cite{F}).
The computational effort is therefore controlled by the process $G$.
We want to consider here schemes $\{\pi^n\}$ which satisfy 
(\ref{sumsq}) together with a mild symmetry requirement:
\begin{equation} \label{sym}\begin{split}
& \mathbb{E}[(\Delta^n_{m+1}W^j)^3 | \mathcal{F}^n_{m}] = 0 ,  \ \ 
\mathbb{E}[L^{n,i,j}_{m+1} | \mathcal{F}^n_{m}] = 0, \ \ 
\mathbb{E}[|\Delta^n_{m+1}W^j|^4 | \mathcal{F}^n_{m}] =
 \frac{1}{n^2} \frac{H^{n}(\pi^n_m)}{G^n(\pi^n_m)}\\
& \textrm{where }L_m^{n,i,j} = \int_{\pi^n_{m-1} }^{\pi^n_{m} }
 (W^i(s) - W^i(\pi^n_{m-1}))(W^j(s) - W^j(\pi^n_{m-1}))\mathrm{d}s
\end{split}
\end{equation}
for each $n$, $m$ and $1 \leq i,j\leq d$ with $i \neq j$, and further
\begin{equation}\label{H}
\begin{split}
&\sup_{ 0 \leq s \leq \tau }|H^{n}(\pi^n(s)) - H(s)| \to 0, \\
&n^2 \sum_{m=1}^{N^n_\tau}\mathbb{E}[|\Delta_{m}\pi^n|^4 |
\mathcal{F}^n_{m-1}] \to 0, \ \ 
n^2 \sum_{m=1}^{N^n_\tau}\mathbb{E}[|\Delta_{m}\pi^n|^6 |
\mathcal{F}^n_{m-1}] \to 0
\end{split}
\end{equation}
in probability as $n\to \infty$,
where $H^{n}$ and $H$ are some progressively measurable processes. Note that the above conditions hold with $G^n = G$ and $H^n = H = 3/G$ if $\pi^n$ is given by (\ref{pred}) with a continuous adapted process $G$.
We can treat a scheme such as
\begin{equation*}
\pi^n_0 = 0, \ \ 
\pi^n_{m + 1} = \pi^n_m + \frac{1}{n (G({\pi^n_m}) \wedge n)}
\end{equation*}
as well by taking $G^n = G \wedge n$.

We now state our central limit theorem for the error process using discretisation schemes as above.
Our standing assumption is that the Euler--Maruyama approximation $X^n$ is well-defined by (\ref{em}), i.e., $X^n(\pi^n_m)$ keeps inside the domain of $f$ up to certain time $t > 0$. Fix such $t>0$. The error process $U^n$ is then well-defined by (\ref{ep}) up to the time $t$.

 \begin{thm} \label{clt}
Let $\{\eta_k\}$ be an increasing sequence of stopping times with
$$\lim_{k \to \infty}\mathbb{P}[\eta_k < t] = 0.$$
Let $\{\pi^n\}$ be a sequence of partitions satisfying, for all $k$, 
(\ref{sumsq}), (\ref{sym}) and (\ref{H}) with locally integrable
  positive processes $G$ and $H$ and  $\tau = \eta_k$. 
 If $\sup_{0 \leq s \leq t} |U^n(s)|$ is tight, then
 the   $C([0,t],\re^p)$-valued random sequence
 $U^n$ converges in law to the solution $U$ of (\ref{SDEU}), where 
$Z$ is given by $Z^{0,j} = Z^{j,0}  = Z^{0,0} = 0$, 
\begin{equation}\label{Z}
 Z^{i,j}(\cdot) =\frac{1}{\sqrt{6}}\int_0^{\cdot} \sqrt{H(s)} \mathrm{d}\hat{W}^{i,j}(s), \ \ 
1 \leq i, j \leq d
\end{equation}
and $\hat{W} = \{\hat{W}^{i,j}\}$ is a $d^2$-dimensional standard Brownian
 motion independent of $W$.
In particular (\ref{Uint}) holds.
\end{thm}

The proof of this result is given in Section~4. 
\begin{remark}\upshape
{The symmetry condition (\ref{sym}) is not essential for a central
 limit theorem to hold. The limit law $Z$ is still a conditionally
 Gaussian semimartingale with anisotropy and drift determined by the asymmetry of the scheme.}
However, we chose to concentrate on the given framework for number of reasons. 
First, asymmetric schemes are in general harder to implement.
Second, whether or not an asymmetric scheme is superior to the standard equidistant
one depends on the coefficients of the SDE (\ref{sde}) which goes against our main purpose of finding  
a generic method which uniformly improves on the standard method.
Third, an asymmetric scheme induces an asymptotic bias and an additional
source of randomness in the limit, which is not preferable from
practical point of view. This bias can be corrected but the correction  makes
implementation more complicated. See Fukasawa~\cite{F}
for related results in the one-dimensional case.
\end{remark}

\begin{remark} \upshape
{ The process $U^n$ we are considering is a normalized
 error process for a continuous version of the
 Euler-Maruyama scheme. 
The whole path generation of $X^n$ requires
a continuous evaluation of a  Brownian motion that is not
 computationally feasible.
Therefore, although Theorem~\ref{clt} ensures the
 convergence of $U^n$ on $C[0,t]$, only relevant in practical applications 
is its finite dimensional convergence
 that follows from Theorem~\ref{clt} as a corollary.
For a simulation of a random variable which depends on the whole path of
 $X$, usually we do linear or piecewise constant interpolation using
 finite evaluations $X^n(\pi^n_m)$, $m=0,1,\dots$, which are
functions of  $\{W(\pi^n_m)\}_m$.
It remains for future research to establish limit
 theorems for those interpolated schemes.
}
\end{remark}

\section{Asymptotically efficient schemes}

We apply now Theorem \ref{clt} to compare and study different discretisation schemes. We establish a lower bound on the expected squared asymptotic error $\E[|U_t|^2]$ and exhibit an efficient scheme which attains the bound. We start however with the usual benchmark given by the classical schemes with Gaussian increments $ \Delta^\pi_mW$.

\subsection{Gaussian schemes}
Consider first $\{\pi^n\}$ defined by (\ref{pred}) with a positive
continuous adapted process $G$.
As already mentioned,
(\ref{sumsq}), (\ref{sym}) and (\ref{H})
are satisfied with $G^n = G$ and
$H^n = H = 3/G$ for any finite stopping time $\tau$.
By the standard theory of the Euler--Maruyama scheme (see e.g.\ Kloeden
and Platen~\cite{KP}) if, say, $G=1$, $\mathbb{D} = \mathbb{R}^p$ and 
the first order derivatives of $f$ are bounded, then
\begin{equation}\label{unifl2}
\sup_{n} \mathbb{E}[ \sup_{0 \leq s \leq t} |U^n(s)|^2] < \infty
\end{equation}
for each $t > 0$, which implies that $\sup_{0 \leq s \leq t} |U^n(s)|$
is tight.
By a localization argument we can easily conclude this tightness without the
restriction $G=1$.
Then we can apply Theorem~\ref{clt} to have
\begin{equation}\label{GaussianLimit}
U(t) = Y(t) \int_0^t Y(s)^{-1} \mathrm{d}F(s), \ \ 
\mathrm{d}F^i(s) =   - \frac{1}{\sqrt{2}}
\sum_{j,l=1}^d \sum_{k=1}^p
\frac{\partial_kf^i_j(X(s))f^k_l(X(s))}{\sqrt{G(s)}} \mathrm{d}\hat{W}^{l,j}(s)
\end{equation}
as the limit of $U^n$,
where $Y = \{Y^{a,b}; 1 \leq a, b \leq p\}$ is the solution of
(\ref{SDEY}).

Next consider $\pi^n_m = g(m/n)$, where $g$ is a continuously
differentiable increasing function with $g(0)=0$.
Since $n \Delta \pi^n_m \sim g^\prime(g^{-1}(\pi^n_{m-1}))$ as $n \to
\infty$, (\ref{sumsq}), (\ref{sym}) and (\ref{H})
are satisfied with $G^n(s) = 1/(g(g^{-1}(s) + 1/n)-s)$, 
$G = 1/g^\prime \circ g^{-1}$, $H^n = 3/G^n$ and $H = 3/G$.
This scheme has essentially the same structure as (\ref{pred}) and 
the tightness of $\sup_{0 \leq s \leq t} |U^n(s)|$ is verified by the
same manner. The limit of $U^n$ is also given by (\ref{GaussianLimit}).

Both of these schemes use conditionally Gaussian increments 
$\Delta^n_m W$. They are relatively easy to implement and widely used in practice.
The particular case $G=1$ is the standard equidistant
scheme. However, when the computational effort is measured by (\ref{effort}), they turn out to be inefficient in terms of
asymptotic error magnitude.

\subsection{Efficient scheme}
By Theorem~\ref{clt}, the asymptotic error distribution depends on
$\{\pi^n\}$ only via $H$ in (\ref{H}), while the computational effort in \eqref{effort} 
does only via $G$. The smaller the $H$, the smaller the error. 
{ In particular, under an additional assumption that $U^n(t)$ is uniformly square
integrable, we have}
\begin{equation}\label{Th}
 \E\left[n|X^n(t)-X(t)|^2\right]=\E\left[|U^n(t)|^2\right]\longrightarrow \E\left[  \int_0^t \Theta(t,s) H(s) \mathrm{d}s\right]
\end{equation}
with $\Theta(t,s)>0$ defined purely in terms of $\{X_s\}_{0\leq s \leq t}$,
where we used (\ref{Uint}) and the independence of $W$ and $\hat W$. 
{ See Remark~\ref{expl} below for an explicit expression of $\Theta.$}
In the following theorem we give a lower bound on $H$, for a given $G$, and exhibit a scheme which attains it. 

\begin{thm}\label{thm:efficient}
Let $t > 0$ and a locally integrable adapted process $G$ be given. 
Let $\{\eta_k\}$ be an increasing sequence of stopping times with
$\lim_{k \to \infty}\mathbb{P}[\eta_k < t] = 0$.
Let $\{\pi^n\}$ be a sequence of partitions satisfying 
 (\ref{sumsq}), (\ref{sym}) and (\ref{H}) 
with a locally integrable process $H$ and  $\tau = \eta_k$ for all $k$. 
Then, 
\begin{equation}\label{lb}
 H(u) \geq \frac{3d}{(d+2)G(u)} \quad \mathrm{d}u\otimes \mathrm{d}\mathbb{P}\text{ a.e. on }[0,t]\times \Omega.
\end{equation}
If $G$ is positive and continuous on $[0,t]$, 
then the sequence $\{\pi^n\}$ defined as
\begin{equation}\label{sd}
\pi^n_0 = 0, \ \ 
\pi^n_{m+1} = \inf\left\{t > \pi^n_m: |W(t) -W(\pi^n_m)|^2
 = \frac{d}{n G(\pi^n_m)}\right\}.
\end{equation}
satisfies  (\ref{sumsq}), (\ref{sym}) and (\ref{H}) 
for $\tau= \hat{\eta}_k$ and attains the equality in (\ref{lb}),
where
\begin{equation*}
 \hat{\eta}_k = \inf\{s >0: G(s) >k \text{ or } G(s) < 1/k\}.
\end{equation*}
\end{thm}
\begin{remark}
\upshape
{ Theorem~\ref{thm:efficient} states a relation between $G$ and $H$
 under
 (\ref{sumsq}), (\ref{sym}) and (\ref{H}). 
To apply Theorem~\ref{clt}, we need additionally to assume
 that the sequence $\sup_{s \in [0,t]}|U^n(s)|$ is tight.
By a standard argument using Gronwall's inequality,
the tightness is verified via the $L^2$ boundedness at least when
$\mathbb{D} =
 \mathbb{R}^p$, $f$ is Lipschitz continuous, 
and 
 (\ref{sumsq}), (\ref{sym}) and (\ref{H}) hold for $\tau = t$
with $H^n(\pi^n_m)/G^n(\pi^n_m)$ bounded uniformly in $n$ and $m$.}
\end{remark}
\begin{remark}
 \upshape
Compared with the standard equidistant scheme,
the error distribution for the scheme (\ref{sd}) uniformly shrinks ;
denoting by $U_{\text{effi}}$ and 
$U_{\text{Gauss}}$ the limits of $U^n$ associated with
(\ref{sd}) and (\ref{pred}) respectively, we conclude
\begin{equation*}
U_{\text{effi}} \stackrel{\textit{law}}{=} \sqrt{\frac{d}{d+2}}U_{\text{Gauss}}.
\end{equation*}
Both (\ref{sd}) and (\ref{pred}) require the same computational effort
(\ref{effort}).
This refines a result for one-dimensional case given in Fukasawa~\cite{F}.
\end{remark}
\begin{remark} \label{expl} \upshape
{ 
An explicit expression of $\Theta$ in (\ref{Th}) can be given as follows.
 By (\ref{Uint}), together with the fact that $\hat{W}$ is independent
 of $W$,
\begin{equation*}
\begin{split}
& \mathbb{E}\left[
|U(t)|^2
\right]=
\mathbb{E}\left[
\sum_{i,j,k,l,m} Y^{i,j}(t)Y^{i,k}(t)\int_0^t
 \bar{Y}_{j,l}(s)\bar{Y}_{k,m}(s)
\mathrm{d}\langle F^l, F^m \rangle(s)
\right]\\
& \mathrm{d}
\langle F^l, F^m \rangle(s) = \frac{1}{6}H(s)\sum_{a,b=1}^d \left(\sum_{c}
\partial_cf^l_a(X(s)) f^c_b(X(s))\right)
\left(\sum_{c}
\partial_cf^m_a(X(s)) f^c_b(X(s))\right) \mathrm{d}s,
 \end{split}
\end{equation*}
where $\bar{Y} = \{\bar{Y}_{a,b}\} =  Y^{-1}$ and so,
\begin{equation*}
\begin{split}
 \Theta(t,s) = \frac{1}{6}
\sum_{i,j,k,l,m} & Y^{i,j}(t)Y^{i,k}(t) \bar{Y}_{j,l}(s)\bar{Y}_{k,m}(s)
\\ &\sum_{a,b=1}^d \left(\sum_{c}
\partial_cf^l_a(X(s)) f^c_b(X(s))\right)
\left(\sum_{c}
\partial_cf^m_a(X(s)) f^c_b(X(s))\right). 
\end{split}
\end{equation*}
It should be noted that $\Theta (t,s) = 0$ if $f^i_j$, $1 \leq j \leq d$
 are constant. In this case a comparison of discretisation schemes has
 to be based on the limit law of order $1/n$ instead of $1/\sqrt{n}$.
Note also that when $d=1$, or more generally, the coefficients are
 ``commutative'', then the Milstein scheme is feasible and achieves
the better rate $1/n$ of convergence;
see Yan~\cite{Yan} and M\"uller-Gronbach~\cite{MG1}.
}
\end{remark}
\begin{remark} \upshape
 We note that the different discretisation schemes we consider mirror different pathwise constructions of the stochastic It\^o integral. The equidistant scheme in \eqref{pred} with $G\equiv 1$ (or other deterministic function) is akin to the approximation of stochastic integral discussed in F\"ollmer~\cite{fol81}. The asymptotically efficient scheme in \eqref{sd} in contrast corresponds to discretising the path along ``Lebesgue type partition", as in Vovk~\cite{vov12} and Perkowski and Pr\"omel~\cite{PePr14}, see also Davis et al.~\cite{DOS} for a discussion and more general partitions.
\end{remark}

To prove \eqref{lb}, it is sufficient to establish the corresponding inequality
for the random variables $H^n(\pi^n_m)$ and $G^n(\pi^n_m)$.  
Given (\ref{sym}), this is equivalent to
$$ \sum_{j=1}^d \mathbb{E}\left[|\Delta^n_{m+1}W^j|^4 | \mathcal{F}^n_{m}\right] \geq \frac{3d^2}{d+2} \left(\mathbb{E}\left[\Delta_{m+1}\pi^n | \mathcal{F}^n_{m}\right]\right)^2,\quad \forall n,m.$$
It follows that Theorem \ref{thm:efficient} is a direct consequence of the following lemma.
\begin{lem} \label{eff}
Let $a > 0$ and 
\begin{equation*}
 Q(v) = \sum_{j=1}^d |W^j(v)|^4.
\end{equation*}
Then 
\begin{equation}\label{min}
 \min \left\{\mathbb{E}[Q(\tau)]:  \tau \text{ is a stopping time with
  } \mathbb{E}[\tau] = a \right\} = \frac{3d^2a^2}{d+2}
\end{equation}
and the minimum is attained by
\begin{equation}\label{sphere}
 \tau = \inf\left\{t > 0: |W(t)|^2 = da\right\}.
\end{equation}
\end{lem}

The proof of the Lemma is deferred to Section~4. 
We discuss now how to implement (\ref{sd}).
Since $|W|$ and $W/|W|$ are independent, conditionally on $\mathcal{F}^n_m$,
$\Delta^n_m W$ is independent of $\Delta_n \pi^n$ and
uniformly distributed on the sphere with radius 
\begin{equation*}
 \sqrt{\frac{d}{nG(\pi^n_{m-1})}}.
\end{equation*}
By a scaling property,
\begin{equation*}
 \Delta^n_m W \sim  \sqrt{\frac{d}{nG(\pi^n_{m-1})}} \frac{N}{|N|},  \ \ 
\Delta_m \pi^n \sim \frac{d}{nG(\pi^n_{m-1})} \tau_1
\end{equation*}
where $N \sim \mathcal{N}(0,I_d)$ and 
$\tau_1$ is defined by (\ref{sphere}) with $a = 1/d$, which
has  the same law as
the hitting time of $1$ for the
$d$-dimensional Bessel process starting from $0$.
Generating $\Delta^n_m W$ is well discussed and for $\Delta_m \pi^n$,
it suffices to develop an efficient algorithm for generating $\tau_1$ by, say, the
acceptance-rejection method.
An explicit form of the distribution function of $\tau_1$  is given by Ciesielski and Taylor~\cite{CT}. The implementation effort and complexity vary with $d$. We do not pursue this further here. 
Instead, we provide an attractive alternative in the next section.

\subsection{Moving sphere scheme}
We adapt here the moving sphere algorithm presented by Deaconu and Herrmann~\cite{D}. The idea is to consider partitions defined by hitting times of a sphere with a radius shrinking in time. The rate at which radius shrinks is adjusted in such a way that both the time step and the spatial increment have explicit distributions which are easy to simulate numerically. Both distributions are non-trivial in the sense that they admit density on some set of positive Lebesgue measure. This is in contrast to the two extreme schemes: the classical equidistant scheme in which time steps are deterministic and the asymptotically efficient scheme of Theorem \ref{thm:efficient} in which the increment $\Delta^\pi_mW$ is concentrated on a (fixed) sphere.

Let $G$ be a positive continuous adapted process on $[0,t]$. We define the sequence of partitions $\{\pi^n\}$ 
via 
\begin{equation}\label{ms}
\pi^n_0 = 0, \ \ 
\pi^n_{m+1} = \inf\left\{s > \pi^n_m; |W(s) -W(\pi^n_m)|^2
 > G^n_m\psi \left(\frac{s -\pi^n_m}{G^n_m}\right)
		  \right\}, 
\end{equation}
where
\begin{equation*}
G^n_m = \frac{1}{n G(\pi^n_m)}, \ \ 
 \psi(v) = dv \log \frac{a}{v}, \ \ a = \left(
1 + \frac{2}{d}
\right)^{1 + d/2}.
\end{equation*}
Since $\psi(a) = 0$, $\Delta_m \pi^n$ is bounded by $aG^n_m$.
Although $\psi(0) = 0$, we can show $\Delta_m \pi^n > 0$ a.s.~by 
the law of iterated
logarithm\footnote{ By Theorem 2.9.23 (i) and (ii) of
Karatzas and Shreve~\cite{KS}, a squared Brownian motion grows at most of
order $t \log \log (1/t)$ for small $t > 0$,
while $\psi$ grows of order $t \log (1/t)$.}.
Since $|W|$ and $W/|W|$ are independent, conditionally on
$\mathcal{F}^n_{m-1}$,
\begin{equation*}
(\Delta_m \pi^n, \Delta^n_m W)
 \sim 
\left(G^n_{m-1}\tau_\psi ,
 \sqrt{G^n_{m-1}\psi (\tau_\psi )} \frac{N}{|N|}
\right),
\end{equation*}
where $N \sim \mathcal{N}(0,I_d)$ and 
\begin{equation*}
 \tau_\psi = \inf\{s > 0; |W(s)|^2 > \psi(s)\}.
\end{equation*}
Now we show that generating a random variable with the same distribution
as $\tau_\psi$ is quite easy.
By Proposition~2 of Deaconu and Herrmann~\cite{D},
the density of $\tau_\psi$ is given by
\begin{equation*}
s \mapsto  \frac{1}{\Gamma(d/2)2^{d/2}a^{d/2}s} \left|
ds \log \frac{a}{s}
\right|^{d/2}.
\end{equation*}
Remark that a method of Chen et al.~\cite{C} can be applied to prove this
with a slight modification.
As shown by Proposition~A.1 of Deaconu and Herrmann~\cite{D},
we have then that 
\begin{equation*}
 \tau_\psi \sim  a e^{-Z},
\end{equation*}
where $Z$ is a random variable which follows the 
 Gamma distribution with shape $1 + d/2$ and scale $2/d$.
Since $|N|^2$ is independent of $N/|N|$ and follows 
the 
 Gamma distribution with shape $d/2$ and scale $2$,
we can use $|N|^2$ to generate $Z$ as
\begin{equation*}
 Z \sim \frac{1}{d}(|N|^2 + 2E),
\end{equation*}
where $E$ is an exponentially distributed random variable with mean 1
which is independent of $N$. Thus we have
\begin{equation*}
(\Delta_m \pi^n, \Delta^n_m W)
 \sim 
\left(G^n_{m-1}a e^{-Z} ,
 \sqrt{G^n_{m-1}ad Ze^{-Z}} \frac{N}{|N|}
\right)
\end{equation*}
conditionally on $\mathcal{F}^n_{m-1}$.

Now we show that (\ref{sumsq}), (\ref{sym}) and (\ref{H}) are
satisfied. 
Observe that
\begin{equation*}
\begin{split}
& \mathbb{E}[\Delta_{m+1} \pi^n|\mathcal{F}^n_{m}] =
  aG^n_m\mathbb{E}[e^{-Z}] = G^n_m = \frac{1}{nG(\pi^n_m)},\\
&\mathbb{E}[|\Delta_{m} \pi^n|^k|\mathcal{F}^n_{m-1}] =
  a^k|G^n_m|^k\mathbb{E}[e^{-kZ}] = O(n^{-k}),\\
&
 \frac{1}{n}\mathbb{E}[N^n_\tau] = \mathbb{E}\left[
\sum_{m=1}^{N^n_{\tau}}
G(\pi^n_{m-1})\mathbb{E}[\Delta_{m} \pi^n|\mathcal{F}^n_{m-1}]
\right] \leq K \mathbb{E}\left[
\sum_{m=1}^{N^n_{\tau}}\Delta_m \pi^n
\right] \leq K(t + K/n),
\end{split}
\end{equation*}
where
\begin{equation*}
 \tau = \eta_K = t \wedge \inf\{s > 0: G(s) > K \text{ or } G(s) < 1/K\}
\end{equation*}
for any $K \in \mathbb{N}$.
It follows then
\begin{equation*}
 n^{k-2} \sum_{m=1}^{N^n_\tau} \mathbb{E}[|\Delta_{m}
  \pi^n|^k|\mathcal{F}^n_{m-1}] \to 0
\end{equation*}
in probability for all $k \geq 2$ and $\tau = \eta_K$.
Further by the same argument as in the proof of Lemma~\ref{eff},
we have that
\begin{equation*}
\begin{split}
\mathbb{E}[|\Delta^n_{m+1}W^j|^4 | \mathcal{F}^n_m] & = 
\frac{1}{d}\mathbb{E}[\sum_{j=1}^d|\Delta^n_{m+1}W^j|^4 |
 \mathcal{F}^n_m]  = \frac{3}{d(d+2)} \mathbb{E}[|\Delta^n_{m+1} W |^4| \mathcal{F}^n_m] \\
& = \frac{3|G^n_m|^2}{d(d+2)} \mathbb{E}[\psi(\tau_\psi)^2]  = \frac{3d|G^n_m|^2}{d+2} a^2 \mathbb{E}[Z^2e^{-2Z}] \\ 
&= \frac{ 3 |G^n_m|^2  (d+2)^{d+2}}{d^{d/2}(d+4)^{2+d/2}} = \frac{1}{n^2} \frac{H(\pi^n_m)}{G(\pi^n_m)},
\end{split}
\end{equation*}
where
\begin{equation}\label{eq:r(d)_ms}
H =  
 \frac{3r(d)}{G}, \ \ 
r(d) = \frac{  (d+2)^{d+2}}{d^{d/2} (d+4)^{(d+4)/2}}.
\end{equation}

\begin{figure}
\begin{center}
 \includegraphics[width=8cm]{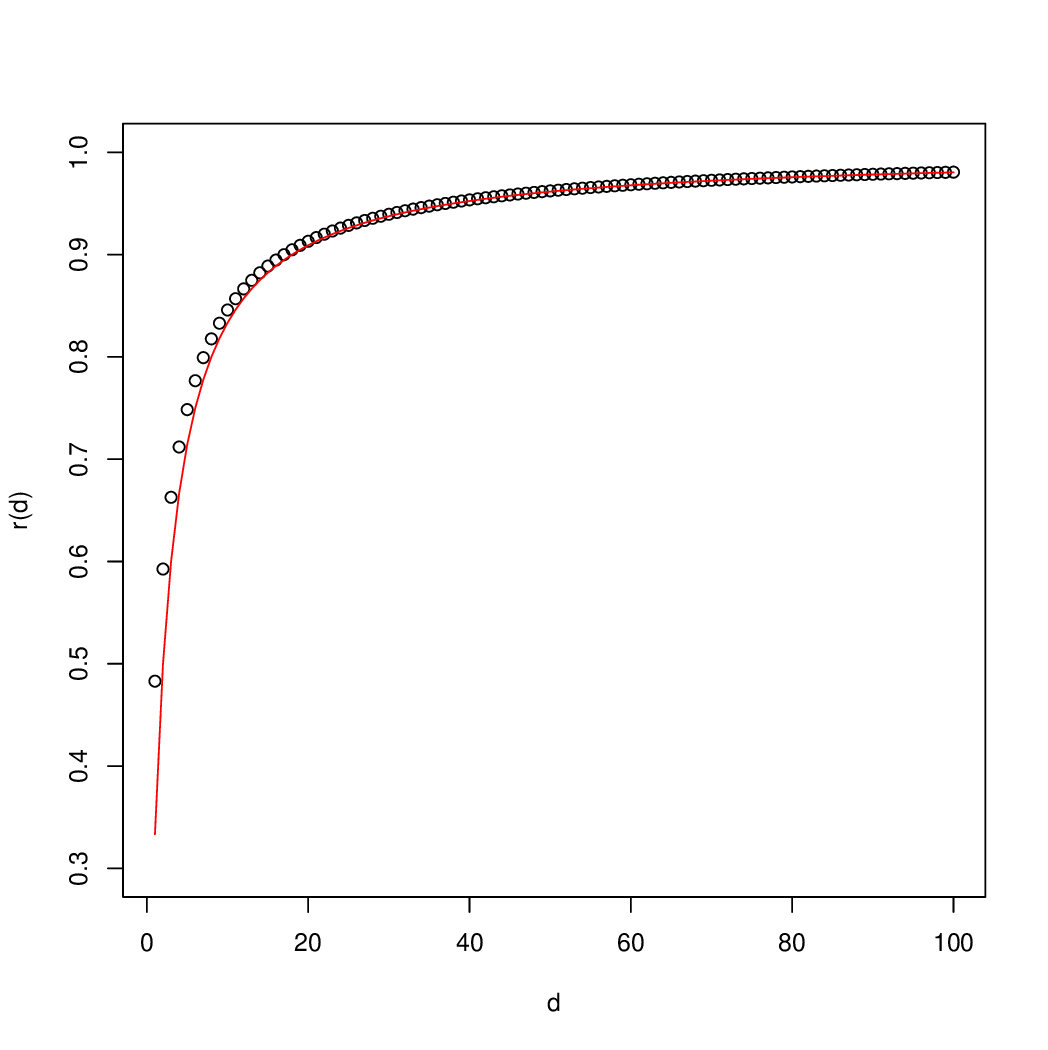}
\end{center}
\caption{Reduction ratio $r(d)$ in \eqref{eq:r(d)_ms} for the moving sphere scheme.}\label{rr}
\end{figure}

Figure~\ref{rr} plots the reduction ratio $r(d)$ with 
efficiency bound $d/(d+2)$ from (\ref{lb}) in red.
As is clearly seen, $r(d) < 1$, which means (\ref{ms}) is superior to
Gaussian schemes. 
We can show\footnote{
Let $f(d) = \log(r(d)(d+1)/d)$. Then,
$$f^\prime(d) = \log (d+2) -\frac{1}{2}\log (d+4) -\frac{1}{2}\log d -
\frac{1}{d(d+1)}, \ \ 
f^{\prime\prime}(d) = \frac{-2d^3 + 5d^2 + 18d +
8}{d^2(d+1)^2(d+2)(d+4)}.$$
It is easy to see $f^\prime(d) \to 0$ as $d \to \infty$ and 
$f^{\prime\prime}(d) < 0$ for $d\geq 5$.
This means $f^\prime (d) > 0$ for $d \geq 5$.
Since $f(d) \to 0$ as $d\to\infty$, this further implies that $f(d) < 0$
for $d \geq 5$. The same inequality can be directly checked for
$d=1,2,3,4$.
Thus we conclude $f(d) < 0$, which is equivalent to $r(d) < d/(d+1)$.  
}
 that $r(d) < d/(d+1)$. 
This inequality implies that the moving sphere
scheme keeps the advantage even after
 taking into account that the proposed method 
requires generating one additional exponential random variable in each step.

Further from Figure~\ref{rr}, we find that 
the reduction ratio $r(d)$ is close to 
the best possible value $d/(d+2)$ attained by (\ref{sd}).
Taking computational costs for generating increments 
$\Delta_m\pi^n$ into account,  (\ref{ms}) is likely to be more efficient than (\ref{sd}) for most applications.

Finally, the moving sphere scheme has an additional advantage that both the
increments $\Delta_m \pi^n$ and $\Delta^n_m W$ are bounded.
More precisely
\begin{equation}
\Delta_m\pi^n \leq aG^n_m\quad \textrm{and}\quad |W(s)-W(\pi_m^n)|\leq \frac{da}{\mathrm{e}}G^n_m
\end{equation}
for all $m$ and $n$.
In consequence, by changing $G$ in an adapted way, we can control the size of each increments of $X^n$. This is  necessary to deal with the SDE on a bounded domain, see Deaconu and Herrmann~\cite{D} or Milstein and Tretyakov~\cite{MT}. Further, this also allows to adapt the scheme to obtain a greater accuracy of certain path information. Consider for example pricing of a barrier option -- we may monitor the barrier crossing of our simulated paths with arbitrary accuracy. This would replace the Brownian bridge correction usually combined with the equidistant scheme. We note that similar consideration motivated recently Gassiat et al.\ \cite{Gassiatetal:2014} who used Root's barrier hitting times for a $1$-dimensional Brownian motion in discretisation schemes.

\subsection{Implementation and numerical experiments}
For implementation purposes, note that the moving sphere scheme, or the asymptotically efficient scheme in \eqref{sd}, have random time steps and we would typically overshoot the time horizon, i.e.\ have that $\pi^n_{N_t^{\pi^n}}>t$. However in the case of moving sphere scheme this is easily corrected. Indeed, the bound $\Delta_m\pi^n \leq aG^n_m$ means we can run the time steps until $t-\pi_n^m\geq aG^n_m$ and when this fails we may either chose to decrease $G_m^n$ or else continue with the equidistant scheme. This would involve at most $6$ equidistant time steps in $d=1$ and at most $4$ for $d\geq 2$. 

We report now a brief numerical study.
We compare the moving sphere scheme with $G=1$
 and the standard equidistant scheme. Consider a two dimensional SDE
\begin{equation*}
 \begin{split}
\mathrm{d}X^1(s) &= 
- \left(\frac{\tan ((X^1(s)+X^2(s))/2) }{(1 + \tan^2((X^1(s) + X^2(s))/2))^2} +
 \frac{\tan ((X^1(s)-X^2(s))/2) }{(1 + \tan^2((X^1(s) - X^2(s))/2))^2} 
\right)\mathrm{d}t \\ & + 
\frac{1 }{1 + \tan^2((X^1(s) + X^2(s))/2)} \mathrm{d}W^1(s) + 
\frac{1 }{1 + \tan^2((X^1(s) - X^2(s))/2)} \mathrm{d}W^2(s)
\\
 \mathrm{d}X^2(s)& =  
- \left(\frac{\tan ((X^1(s)+X^2(s))/2) }{(1 + \tan^2((X^1(s) + X^2(s))/2))^2} -
 \frac{\tan ((X^1(s)-X^2(s))/2) }{(1 + \tan^2((X^1(s) - X^2(s))/2))^2} 
\right)\mathrm{d}t \\ & + 
\frac{1 }{1 + \tan^2((X^1(s) + X^2(s))/2)} \mathrm{d}W^1(s) -
\frac{1 }{1 + \tan^2((X^1(s) - X^2(s))/2)} \mathrm{d}W^2(s)
\end{split}
\end{equation*}
with $X^1(0) = X^2(0) = 0$.
{ It is straightforward to check by It\^o's formula that 
\begin{equation*}
 X^1(s) = \arctan(W^1(s)) + \arctan(W^2(s)), \ \ 
 X^2(s) = \arctan(W^1(s)) - \arctan(W^2(s))
\end{equation*}
satisfies the SDE. By comparing with this explicit solution,
 we can assess the quality of the Euler-Maruyama approximations.}
We approximate the first four moments of the discretisation errors $E^{n,1} := X^{n,1}(1)-X^1(1)$ and $E^{n,2}:= X^{n,2}(1)-X^2(1)$ 
by the Monte Carlo with 1,000,000 paths for both of the schemes.
Normal random variables are generated by the Box-M$\ddot{\text{u}}$ller algorithm and an exponential variable is by $ - \log (U)$, where $U$ is a uniform random variable. Uniform random variables are generated by the Mersenne Twister algorithm.
The Apple Mac mini computer with 2.6 GHz Intel Core i7 took 1 minute and 28 seconds for the equidistant scheme with $n = 625$. For the moving sphere scheme we took $n=435$ which leads to an equivalent  computational time of 1 minute and 26 seconds\footnote{As described above, we follow the moving scheme algorithm until $t-\pi_n^m\geq aG^n_m$ and then we finish with $4$ equidistant steps. This led to an average of $435.9$ steps.}.
These two are therefore almost equivalent in terms of computation time.
Table~\ref{m1}
and \ref{m2} 
 report the first four moments of $E^{n,1}$ and $E^{n,2}$ respectively.

\begin{table}[htbp]
 \begin{center}
 \begin{tabular}{|c|c|c|c|c|c|}
\hline
 & $\E[E^{n,1}]$ & $\E[|E^{n,1}|^2]$ & $\E[(E^{n,1})^3]$ & $\E[|E^{n,1}|^4]$\\
\hline
equidistant ($n=625$) &  8.1 $\times 10^{-6}$ & 0.00033 & 1.1 $\times 10^{-8}$ & 4.1 $\times 10^{-7}$ \\
\hline
moving sphere ($n=435$) 
&  -4.7 $\times 10^{-6}$ & 0.00028 & 1.9 $\times 10^{-8}$ & 2.9 $\times 10^{-7}$ \\
\hline
 \end{tabular} 
 \end{center}
\caption{$E^{n,1} = X^{n,1}(1) - X^1(1)$}\label{m1}
\end{table}

\begin{table}[htbp]
 \begin{center}
 \begin{tabular}{|c|c|c|c|c|c|}
\hline
 & $\E[E^{n,2}]$ & $\E[|E^{n,2}|^2]$ & $\E[(E^{n,2})^3]$ & $\E[|E^{n,2}|^4]$\\
\hline
equidistant ($n=625$) &  -1.1 $\times 10^{-5}$ & 0.00033 & -3.3 $\times 10^{-8}$ & 4.1 $\times 10^{-7}$ \\
\hline
moving sphere ($n=435$) 
&  -1.1 $\times 10^{-5}$ & 0.00028 & -7.8 $\times 10^{-10}$ & 2.9 $\times 10^{-7}$ \\
\hline
 \end{tabular} 
 \end{center}
\caption{$E^{n,2} = X^{n,2}(1) - X^2(1)$}\label{m2}
\end{table}
From Table~\ref{m1} and \ref{m2}, we confirm that the moving
sphere scheme provides a better accuracy without increasing computation time.

\section{Proofs}
\subsection{Proof of Theorem \ref{clt}}
Let 
$\epsilon_m = \inf\{|x-y|;x \in \mathbb{K}_m, 
y \in \mathbb{K}_{m+1}^c\} > 0$,
$\hat{\tau}_m = \tau_m \wedge \eta_m$
 and 
\begin{equation*}
 \sigma^n_m = \eta_m \wedge 
\inf\{ u > 0; X(u) \not\in \mathbb{K}_m \text{ or } X^n(u)
  \not\in \mathbb{K}_{m+1}\}.
\end{equation*}
Then
\begin{equation*}
\mathbb{P}[\hat{\tau}_m < t] \leq
 \mathbb{P}[\sigma^n_m < t] \leq 
\mathbb{P}[\hat{\tau}_m < t]
+ \mathbb{P}[\sup_{0 \leq s \leq t}|U^n(s)| \geq \sqrt{n}\epsilon_m].
\end{equation*}
Since $\sup_{0 \leq s \leq t} |U^n(s)|$ is tight by the assumption,
for any $\epsilon > 0$, there exists $K > 0$ such that
\begin{equation*}
 \limsup_{n\to\infty}
\mathbb{P}[\sup_{0 \leq s \leq t}|U^n(s)| \geq K] < \epsilon.
\end{equation*}
It follows that $\sigma^n_m \to \hat{\tau}_m$, 
$1_{\{\sigma^n_m < t\}} \to 1_{\{\hat{\tau}_m < t\}}$ 
in probability and
\begin{equation*}
\lim_{m\to\infty}
\limsup_{n\to \infty} \mathbb{P}[\sigma^n_m < t] = 0.
\end{equation*}
For any continuous bounded function $\varphi$ on $C([0,t],\re^p)$,
\begin{equation*}
 \mathbb{E}[\varphi(U^n)] = 
\mathbb{E}[\varphi(U^n(\cdot \wedge \sigma^n_m))1_{\{\sigma^n_m \geq t\}}] + 
\mathbb{E}[\varphi(U^n)1_{\{\sigma^n_m < t\}}].
\end{equation*}
Therefore it suffices to show
\begin{equation*}
 \mathbb{E}[\varphi(U(\cdot \wedge \hat{\tau}_m))1_{\{\hat{\tau}_m \geq
  t\}}] 
= \lim_{n \to \infty}
 \mathbb{E}[\varphi(U^n(\cdot \wedge \sigma^n_m))1_{\{\sigma^n_m \geq
  t\}}].
\end{equation*}
The coefficient $f$ and its first derivatives
are bounded and uniformly continuous on the compact sets $\mathbb{K}_{m+1}$.
Therefore, for $v \leq \sigma^n_m \wedge t$,
\begin{equation*}
\begin{split}
U^{n,i}(v) =& \sqrt{n} \sum_{j=0}^d
\int_0^v \left\{f^i_j(X^n(\pi^n(s)))-f^i_j(X(s))\right\}\mathrm{d}W^j(s) \\
=& \sqrt{n}  \sum_{j=0}^d 
\int_0^v \left\{f^i_j(X^n(\pi^n(s))) - f^i_j(X^n(s)) + f^i_j(X^n(s))
-f^i_j(X(s))\right\}\mathrm{d}W^j(s) \\
=& - \sqrt{n}\sum_{j,k} \int_0^v  \partial_k f^i_j(X^n(\pi^n(s)))
(X^{n,k}(s) - X^{n,k}(\pi^n(s))) \mathrm{d}W^j(s)
 \\ & + \sqrt{n}\sum_{j,k}
\int_0^v  \partial_k f^i_j(X(s))(X^{n,k}(s)-X^k(s))
 \mathrm{d}W^j(s) + o_p(1)\\
=& - \sqrt{n}\sum_{j,k,l} \sum_{m=0}^{\infty} \partial_k f^i_j(X^n(\pi^n_m))
 f^k_l(X^n(\pi^n_m))
\int_{\pi^n_m \wedge v}^{\pi^n_{m+1} \wedge v}
 (W^l(s) - W^l(\pi^n_m))\mathrm{d}W^j(s) \\ &+ \sum_{j,k}
\int_0^v \partial_k f^i_j(X(s))U^{n,k}(s)\mathrm{d}W^j(s) + o_p(1).
\end{split}
\end{equation*}
Denote by $V^{n,i}$ the first of the two terms in the final expression above. 
Put $X^n_m = X^n(\pi^n_m)$ and
\begin{equation*}
 L_m^{n,b,c}(v) = \int_{\pi^n_{m-1} \wedge v}^{\pi^n_{m} \wedge v}
 (W^b(s) - W^b(\pi^n_{m-1}))(W^c(s) - W^c(\pi^n_{m-1}))\mathrm{d}s
\end{equation*}
for $0 \leq b,c\leq d$.
Then $V^n = (V^{n,1},\dots,V^{n,p})$ is
 a  continuous semimartingale with 
quadratic covariation $\langle V^{n,i}, V^{n,j} \rangle_v$
 given by
\begin{equation*}
n \sum_{m=0}^{\infty}
\sum_{a=1}^d\sum_{b,c=0}^d
\sum_{k,l=1}^p\partial_k f^i_a(X^n_m)
\partial_l f^j_a(X^n_m)
f^k_b(X^n_m)
f^l_c(X^n_m) L_{m+1}^{n,b,c}(v)
\end{equation*}
and
\begin{equation*}
\langle V^{n,i},W^j \rangle_v = -
\sqrt{n} \sum_{m=0}^{\infty} \sum_{k=1}^p \sum_{l=0}^d\partial_k f^i_j(X^n_m)
 f^k_l(X^n_m)
\int_{\pi^n_m \wedge v}^{\pi^n_{m+1} \wedge v}
 (W^l(s) - W^l(\pi^n_m))\mathrm{d}s.
\end{equation*}
By Theorem~IX.7.3 of Jacod and Shiryaev~\cite{JS}, if there exists a continuous process $A =
\{A^{i,j}\}$ such that
\begin{equation}\label{eq:key_cnv}
\langle V^{n,i},V^{n,j} \rangle \to A^{i,j}, \ \ 
\langle V^{n,i}, W^j \rangle \to 0
\end{equation}
in probability as $n \to \infty$ for all $i,j$,
then $V^n$ converges $\mathcal{F}$-stably in law 
to a conditionally Gaussian martingale 
$V = (V^1,\dots, V^p)$ 
with $\langle V^i, V^j \rangle = A^{i,j}$.

We will argue below that \eqref{eq:key_cnv} holds and that, 
using $Z$ defined by (\ref{Z}), the limit $V$ is written as
\begin{equation}\label{V}
V^i(v) = - \sum_{a,b=1}^d\sum_{k=1}^p\int_0^v \partial_k f^i_a(X(s))f^k_b(X(s)) \mathrm{d}Z^{b,a}(s).
\end{equation}
The convergence of $V^n$ implies tightness of
$U^n_{\cdot \wedge \sigma^n_m}$ in $C[0,t]$ by Theorem~VI.4.18 of Jacod and Shiryaev~\cite{JS}.
So any subsequence has a further subsequence which converges in law. 
Further it follows from  (\ref{V}) that  the limit of the subsequence is uniquely determined by the SDE 
(\ref{SDEU}). Therefore $U^n_{\cdot \wedge \sigma^n_m}$ 
itself must converges to $U_{\cdot \wedge \hat{\tau}_m}$ stably and we easily conclude.

It remains to establish \eqref{eq:key_cnv}. We do this in two steps.\smallskip\\
{\bf Step 1)}: We first show $\langle V^{n,i}, W^j \rangle \to 0$.\\
By It$\hat{\text{o}}$'s formula
\begin{equation*}
\int_{\pi^n_{m-1}}^{\pi^n_{m}}
 (W^l(s) - W^l(\pi^n_{m-1})) \mathrm{d}s
 =  
\frac{1}{3} (\Delta^n_{m}W^l)^3 -
\int_{\pi^n_{m-1}}^{\pi^n_{m}}
 (W^l(s) - W^l(\pi^n_{m-1}))^2 \mathrm{d}W^l(s)
\end{equation*}
for all $1 \leq l \leq d$. 
The conditional expectations of both terms in the right hand side are $0$ by (\ref{sym}).
Further by (\ref{sumsq}) and (\ref{H}),
\begin{equation}\label{sumtri}
n \sum_{m=1}^{N^n_v} \mathbb{E}[|\Delta_{m}
 \pi^n|^3|\mathcal{F}^n_{m-1}] \leq 
\sqrt{n^2 \sum_{m=1}^{N^n_v} \mathbb{E}[|\Delta_{m}
 \pi^n|^4|\mathcal{F}^n_{m-1}]}
\sqrt{\sum_{m=1}^{N^n_v} \mathbb{E}[|\Delta_{m}
 \pi^n|^2|\mathcal{F}^n_{m-1}]} \to 0
\end{equation}
in probability and
\begin{equation*}
\begin{split}
& \mathbb{E}[|\Delta^n_m W|^6| \mathcal{F}^n_{m-1}]
\leq C
\mathbb{E}[|\Delta_m \pi^n|^3 | \mathcal{F}^n_{m-1}], \\
&  \mathbb{E}[
\int_{\pi^n_{m-1}}^{\pi^n_{m}} (W^l(s)-W^l(\pi^n_{m-1}))^4 \mathrm{d}s | \mathcal{F}^n_{m-1}
] \\
& \leq \frac{1}{15} \liminf_{u \to \infty}
\mathbb{E}[|W^l(u\wedge \pi^n_m)-W^l(\pi^n_{m-1})|^6| \mathcal{F}^n_{m-1}]
\leq C
\mathbb{E}[|\Delta_m \pi^n|^3 | \mathcal{F}^n_{m-1}] 
\end{split}
\end{equation*}
for a constant $C > 0$.
Then by Lemma~A.2 of Fukasawa~\cite{F},
we obtain
\begin{equation*}
 \sqrt{n} \sum_{m=0}^{\infty} \sum_{k=1}^p \partial_k f^i_j(X^n_m)
 f^k_l(X^n_m)
\int_{\pi^n_m \wedge v}^{\pi^n_{m+1} \wedge v}
 (W^l(s) - W^l(\pi^n_m))\mathrm{d}s \to 0
\end{equation*}
in probability for $1 \leq l \leq d$.
To treat the case $l=0$,
observe that
\begin{equation*}
\sqrt{n} \sum_{m=1}^{N^n_v} \mathbb{E}[|\Delta_{m}
 \pi^n|^2|\mathcal{F}^n_{m-1}] \leq 
\sqrt{nN^n_v} \sqrt{\sum_{m=1}^{N^n_v}\mathbb{E}[
|\Delta_{m}
 \pi^n|^4|\mathcal{F}^n_{m-1}]} \to 0
\end{equation*}
 by (\ref{effort}) and (\ref{H}).
It follows then that
$\langle V^{n,i}, W^j \rangle \to 0$
for all $i,j$ again with the aid of Lemma~A.2 of Fukasawa~\cite{F}.
\smallskip\\
{\bf Step 2)}: We show that $\langle V^{n,i},V^{n,j} \rangle$ converges and compute the limit $A^{i,j}$.\\
By It$\hat{\text{o}}$'s formula, (\ref{sym}), we get
\begin{equation*}
\mathbb{E}[L_{m+1}^{n,b,c}(\pi^n_{m+1})| \mathcal{F}^n_m] = 
\mathbb{E}[L^{n,b,c}_{m+1}|\mathcal{F}^n_m] = \frac{\delta^{b,c}}{6n}H^n(\pi^n_m)
 \mathbb{E}[\pi^n_{m+1}-\pi^n_m| 
\mathcal{F}_{\pi^n_m}]
\end{equation*}
for $1 \leq b,c \leq d$,
where $\delta^{b,c}$ is Kronecker's delta.
The terms with $b =0$ or $c=0$ are negligible since
\begin{equation*}
| L^{n,b,c}_{m+1}(\pi^n_{m+1})|\leq
\sqrt{
L^{n,b,b}_{m+1}(\pi^n_{m+1}) 
} 
\sqrt{
 L^{n,c,c}_{m+1}(\pi^n_{m+1}) 
}
\end{equation*}
and 
\begin{equation*}
n \sum_{m=1}^{N^n_v} |\Delta_m \pi^n|^3 \to 0
\end{equation*}
in probability, which follows from (\ref{H}) and (\ref{sumtri}) 
by using Lemma~A.2 of Fukasawa~\cite{F}.
Therefore,
\begin{equation*}
\langle V^{n,i},V^{n,j} \rangle_v
\to \frac{1}{6}\sum_{a,b=1}^d \sum_{k,l=1}^p    \int_0^v
\partial_k f^i_a(X(s))
\partial_l f^j_a(X(s))
f^k_b(X(s))
f^l_b(X(s)) H(s)\mathrm{d}s
\end{equation*}
again by Lemma~A.2 of Fukasawa~\cite{F}. This completes the proof.

\subsection{Proof of Lemma \ref{eff}}
Let $\tau$ be a stopping time with $\mathbb{E}[\tau] = a$.
Then $W^j_{\cdot \wedge \tau}$ are uniformly integrable martingales and
so, for any $t > 0$ by Jensen's inequality,
\begin{equation*}
 \mathbb{E}[Q(\tau)|\mathcal{F}_{\tau \wedge t}] 
\geq \sum_{j=1}^d
|\mathbb{E}[W^j(\tau)| \mathcal{F}_{\tau \wedge t}]|^4 =
Q(\tau \wedge t).
\end{equation*}
Therefore for (\ref{min}), it suffices to show 
\begin{equation*}
 \mathbb{E}[Q(\tau)] \geq \frac{3d^2}{d+2} \mathbb{E}[\tau]^2
\end{equation*}
when $\tau$ is a bounded stopping time.
Let
\begin{equation*}
 S(v) = \sum_{j=1}^d |W^j(v)|^2.
\end{equation*}
Then
\begin{equation*}
 \mathrm{d}Q(v) = 4 \sum_{j=1}^d (W^j(v))^3\mathrm{d}W^j(v) + 
6 S(v)\mathrm{d}v
\end{equation*}
and
\begin{equation*}
 \mathrm{d}S^2(v) = 2S(v)\mathrm{d}S(v) + 4 S(v)\mathrm{d}v = 
4S(v)\sum_{j=1}^dW^j(v)\mathrm{d}W^j(v) + 2S(v)(2+d)\mathrm{d}v.
\end{equation*}
It follows then
\begin{equation*}
 \mathbb{E}[Q(\tau)] = 6 \mathbb{E}[\int_0^\tau S(v)\mathrm{d}v] = 
\frac{3}{2+d}\mathbb{E}[S(\tau)^2]
\geq \frac{3}{2+d}\mathbb{E}[S(\tau)]^2 = 
\frac{3d^2}{2+d}\mathbb{E}[\tau]^2.
\end{equation*}
The equality is attained if and only if $S(\tau)$ is a constant, or
equivalently, $\tau$ is given by (\ref{sphere}). This completes the
proof.
\\

\noindent
{\bf Acknowledgement.} The authors thank the two anonymous reviewers for
their helpful comments and suggestions.

\end{document}